\documentclass[10pt]{article}
\usepackage{amsmath}
\usepackage{amsfonts}
\usepackage{amssymb}
\usepackage[T1]{fontenc}
\usepackage[utf8]{inputenc}
\usepackage{graphicx}
\usepackage{amsthm}
\usepackage{hyperref}
\usepackage{xcolor}
\usepackage{accents}
\usepackage{graphics,wrapfig, graphicx, mathrsfs,xcolor}
\usepackage{mathtools}
\usepackage{amssymb,amscd,graphicx,xcolor,subfig,tikz}
\usepackage{graphics}
\usepackage{indentfirst}
\usepackage{bm, enumerate,xcolor}
\usepackage[dvips]{epsfig}
\usepackage{latexsym}
\usepackage{float}
\usepackage{amsmath}
\usepackage[numbers,sort&compress,square]{natbib}

\usepackage[english]{babel}
\usepackage [autostyle, english = american]{csquotes}
\MakeOuterQuote{"}

\numberwithin{equation}{section}
\theoremstyle{plain}
\newtheorem*{theorem*}{Theorem}
\newtheorem*{lemma*}{Lemma}

\newtheorem{lemma}{Lemma}[section]

\newtheorem{proposition}[lemma]{Proposition}

\theoremstyle{definition}

\newtheorem{remark}[lemma]{Remark}

\newtheorem*{example*}{Example}
\newtheorem*{er*}{Examples and Remarks}

\usepackage{geometry}
\newcommand{\dth}{\partial_{\theta}}

\newcommand{\dr}{\partial_{r}} 
 
\newcommand{\dt}{\partial_{t}}

\newcommand{\til}{\tilde}

\begin{document}

\title{BMO Vorticity Strong Ill-Posedness for the 2d Incompressible  Euler Equation}
\author{Karim R. Shikh Khalil}
\maketitle
\abstract{For the  two dimensional incompressible Euler equation, Yudovich, and independently Serfati, extended the classical existence and uniqueness class to include vorticity with mild singularity. This class is sometimes called the generalized Yudovich class. Then, Vishik extended the uniqueness class to include vorticity in BMO, allowing for stronger singularities than the generalized Yudovich class.  In this  work, we prove  strong ill-posedness of vorticity in BMO. }

\section{Introduction}
The two dimensional incompressible Euler equation  in vorticity form is: 
 \begin{equation}\label{EulerVorticty}
\begin{split}
\partial_t \omega &+  u\cdot\nabla \omega= 0, \\ 
u&=\nabla^\perp \Delta^{-1} \omega,
\end{split}
\end{equation}  
where $u$ is  the velocity field, and $\omega$ is the vorticity. The classical work of Yudovich \cite{Y1}  proved global \textit{existence} \textit{and} \textit{uniqueness} of integral and bounded vorticity $\omega_{0} \in L^{\infty} \cap L^{1}$. Yudovich \cite{Y2}, and independently Serfati \cite{Se}, extended the global \textit{existence} \textit{and} \textit{uniqueness} class to include mildly unbounded vorticity which is sometimes called the generalized Yudovich class. This generalized Yudovich class includes vorticities with a  $\log(|\log(\cdot)|)$ type singularity. Then, Vishik \cite{MV2} extended the \textit{uniqueness} class of the 2d Euler equation to include vorticities in the Bounded Mean Oscillation ($BMO$) space, thus enlarging  the uniqueness class for 2d Euler since functions in $BMO$ can have a $|\log(\cdot)|$ type singularity. The point to emphasize here is Vishik's result \cite{MV2} gave only \textit{uniqueness}  in $BMO$, so the goal of this work is to prove  strong ill-posedness of vorticity  in $BMO$. Before we get into more details about this work (see subsection \ref{BMOMainResult}), we will give a brief review of the well/ill-posedness works on the 2d Euler equation. The literature on the 2d Euler equation is vast so the review here is not exhaustive.  
 \subsection{Brief overview of well/ill-posedness works}

The well-posedness of the Euler equation was proved in H\"older spaces $C^{k,\alpha}$ by Lichtenstein~\cite{LL} and Gunther~\cite{G}.  Wolibner~\cite{Wo} and H\"older~\cite{HO} proved global well-posedness in two dimensions. Then,   Ebin and Marsden in \cite{EbM} proved  well-posedness in the  Sobolev spaces $H^s$ with $s > \frac{d}{2} + 1$, for compact domains with boundary. Bourguignon and Brezis~\cite{BB} extended the well-posedness to $W^{s,p}$ spaces for $s > \frac{d}{p} + 1$. Further,  Kato in \cite{Ka2} and Kato and Ponce in \cite{KaP}  extended   the  Sobolev spaces well-posedness to $\mathbb{R}^d$.  Vishik in \cite{MV1} obtained global well-posedness in the critical Besov spaces $B^{2/p+1}_{p,1}$ for $1 < p < \infty$. Then, Chae~\cite{C}  extended this to $B^{d/p+1}_{p,1}$ spaces. Pak and Park~\cite{PP1,PP2}, based on Vishik~\cite{MV1}, proved well-posedness for $B^1_{\infty,1}$ and $B^{d+1}_{1,1}$ Besov spaces. Recently, Cozzi and Harrison in \cite{CH} proved local well-posedness in two dimensions for vorticity in the critical  $W^{2,1}$ space. Then, Danchin  and Jeong in \cite{DJ}  extended this to global well-posedness, and proved, in three dimensions, global well-posedness of axisymmetric flows without swirl in $W^{3,1}$. See \cite{CH,DJ} for more discussion.

In all the above works, the gradient of the velocity field is bounded which is the key to propagating the regularity. When the vorticity has lower regularity, for example when it is only bounded, the velocity field is not necessarily Lipschitz. Yudovich, in his  classical work~\cite{Y1}, proved global existence and uniqueness for initial vorticity $\omega_0 \in L^{\infty} \cap L^{1}$. Later on, Yudovich~\cite{Y2}, and independently Serfati~\cite{Se}, proved existence and uniqueness for unbounded vorticity with mild singularities (The generalized Yudovich class). Then, Vishik in \cite{MV2}, as we mentioned, further extended the uniqueness class to a borderline Besov space that includes $BMO$. Bernicot and Keraani in \cite{BK} proved existence and uniqueness in a space strictly between $L^{\infty}$ and $BMO$, called $LBMO$, where the $L$ stands for logarithmic. It is interesting to compare our paper with \cite{BK}. The $LBMO$ space puts more constraints on unbounded functions than $BMO$ so this allowed for the propagation of  regularity. There are also other works constructing special types of solutions for unbounded vorticities. Elling in \cite{Ell1,Ell2} constructed self-similar algebraic spiral solutions. See the works of Shao, Wei, and Zhang in \cite{SWZ} and Choi in \cite{Ch}  for extensions. In addition, there is  work of Elgindi, Murray, and Said~\cite{EMS} which  gave well-posedness class for a certain type of algebraic singularities.

We now discuss the works on strong ill-posedness (see similar review in \cite{SK}). DiPerna and Lions \cite{DL}, using the $2 \frac{1}{2}$-dimensional shear flows, proved strong ill-posedness in the supercritical $W^{1,p}(\mathbb{T}^3)$ spaces for $1 \leq p < \infty$.  Using these shear flows, Bardos and Titi in \cite{BT} proved strong ill-posedness in the H\"older $C^{\alpha}$ space. Further,  Misio\symbol{170}ek and Yoneda~\cite{MY} also used   these flows  to prove strong ill-posedness in the log-Lipschitz type spaces $LL^{\alpha}$ for $0 < \alpha \leq 1$. In the critical spaces, Bourgain and Li~\cite{BL1} in dimensions $d=2, 3$  proved strong ill-posedness in the Sobolev spaces $W^{\frac{d}{p}+1,p}$ for any $1 < p < \infty$, and for Besov spaces $B^{\frac{d}{p}+1}_{p,q}$ for any $1 < p < \infty$ and $1 < q \leq \infty$. Strong ill-posedness in the $C^k$  integer spaces was  proved independently by Elgindi and Masmoudi~\cite{EM} and Bourgain and Li~ \cite{BL2}. Elgindi and Jeong  \cite{EJ1}, later, gave a simplified proof of ill-posedness in the critical Sobolev spaces, in two dimensions. For supercritical regularity in two dimensions, Jeong~\cite{J} showed that a unique solution in $W^{1,p}$ for $1 \leq p < 2$ can continuously loss $W^{1,q(t)}$ regularity where exponent $q(t) < p$ is decreasing. A gap loss of supercritical regularity was proved by C\'ordoba, Mart\'inez-Zoroa, and O\.za\'nski in \cite{CMO}, where they constructed unique solutions in $H^{\beta}$ for $1 < \beta < 2$ that instantly exit $H^{\beta'}$ for all $\beta' > 1+\frac{(3-\beta)(\beta-1)}{2-(\beta-1)^2}$. In three dimensions, Luo~\cite{L} proved strong ill-posedness in the supercritical Sobolev spaces $H^s$ with $0 < s < \frac{5}{2}$, for the velocity field. Furthermore, for vorticity in $H^s$ with $s \in (0, \frac{3}{2})$, Jeong, Mart\'{i}nez-Zoroa, and O\.{z}a\'{n}ski~\cite{JMO} proved   instantaneous and continuous loss of Sobolev regularity. Recently, in two dimensions,   we proved  \cite{SK} strong ill-posedness of the velocity field in the critical Besov spaces $B^{1}_{\infty,q}$ for $1 < q < \infty$. In the setting of Lorentz spaces, strong ill-posedness was proved by Kim and Jeong in \cite{KJ} and recently extended by Bang and Cheskidov in \cite{BC}.

\subsection{Main result}\label{BMOMainResult} 
$BMO$ regularity, in addition to its importance in analysis and other branches of mathematics (see John \cite{John}, John and Nirenberg \cite{JN}, and Fefferman and Stein \cite{FS}), has significant role for the 2d Euler equation. When the vorticity is in $BMO$, it is at the threshold regularity for which solutions of the 2d Euler equation are guaranteed to be unique. Below $BMO$ uniqueness could fail. However, existence is not known. More precisely, Vishik in \cite{MV2}  proved that \textit{if} one has a solution that stays bounded in the $BMO$, then the solution is unique. The goal of this work is to prove strong ill-posedness of the 2d Euler equation in $BMO$. Namely, one can start with arbitrarily small data in $BMO$, and then the solution becomes arbitrarily large in $BMO$ in  arbitrarily small time. Hence, the solution map  is unbounded in  $BMO$.

One of the main challenges in obtaining the strong ill-posedness in $BMO$ is that one needs to work with unbounded vorticities. This is because $L^{\infty}$ is contained in $BMO$, and the $L^{\infty}$ norm is a conserved for the 2d Euler equation. Hence, it is not possible to obtain any norm inflation in $BMO$ using only bounded data, which is different from previous works of strong ill-posedness. Therefore, we need to work with unbounded vorticity which is generally subtle due to the existence and uniqueness issues discussed earlier.   

Our result is obtained by constructing unique compactly supported solutions beyond the Yudovich class. We approximate unbounded  functions in $BMO$ using H\"older-type ones in which we couple the regularity exponent with appropriate rescaling. This approach has the advantage that now we can study unbounded solutions using the H\"older $C^{\alpha}$  class for which we have local (and global) well-posendess. In addition, using this approach, we make new observation about the $BMO$ space and  obtain quantitive estimates in terms of H\"older regularity exponent. The ill-posedness mechanism in $BMO$ is inspired by our ill-posedness  mechanism in the critical Besov spaces \cite{SK}.  We observe that  at the level of vorticity,  the $BMO$ space and the Besov $B^{0}_{\infty,q}$ space for $1<q<\infty$ share structural  similarities. See Subsection \ref{SubSect_Sketch} for more details.
\vskip 0.2 cm
\noindent \textbf{Main Result Statement:}
\begin{theorem*} (Strong ill-posedness for vorticity in $BMO$):  For any $\epsilon>0$, and $0<\alpha<1$, there exist compactly supported initial vorticity  $\omega_0^{\alpha,\epsilon} $ on $\mathbb{R}^2$  such that the unique solution to the 2d Euler equation  satisfies the following: 
$$|\omega_0|_{BMO}=\epsilon, \, \text{and} \,  \sup_{0\leq t \leq T(\alpha)}|\omega(t)|_{BMO} \geq   \log(c\log(c |\log(\alpha)|)), $$ 
  where 
  $T(\alpha)=c \alpha \log(c |\log(\alpha)|) $ and c is independent of $\alpha$.  
  \end{theorem*}

  \begin{remark} Observe that  as $\alpha \rightarrow 0$, $  |\omega(t)|_{BMO} \rightarrow \infty$  and $T(\alpha)= c \alpha \log(c |\log(\alpha)|)   \rightarrow 0$. Therefore,  we have the  strong ill-posedness result in $BMO$. 
   \end{remark}

\subsection{Sketch of the Proof}\label{SubSect_Sketch}

 There are  three main ideas in our result: 
 
  \textbullet \textit{\, Approximating  $BMO$ functions  using  H\"older  type functions and identifying two different functions in $BMO$:} Here we approximate functions in the $BMO$ space by H\"older-type ones where  we couple the H\"older regularity exponent  with appropriately chosen rescaling. This allows us to study unbounded functions using H\"older-type one, for which the well-posedness theory is known and easier to work with.  In particular, we can obtain quantitive estimate in terms of H\"older exponent. For example, for any $0<\alpha<1$, we consider  following two class of  functions:   
  \begin{equation}\label{BMOFuns}
f(r)= \frac1{\beta(\alpha)}r^{\alpha}, \quad \text{and}   \quad \,  g(r)= \frac{1}{\beta(\alpha)}r^{\alpha} \cos(\beta(\alpha)\log(r)), 
\end{equation}
compactly supported on the ball of radius one. The factor $\beta(\alpha)$ is any continuous function of $\alpha$ satisfying $\lim_{\alpha \rightarrow 0} \beta(\alpha)=0$, and $\frac{\alpha}{\beta(\alpha)}=O(1)$. For instance, $\beta(\alpha)=\alpha \,\,  \text{or}  \,\,   \beta(\alpha)=|\log(\alpha)|^{-1}$ would work. The functions in \eqref{BMOFuns} are $O(1)$ in  $BMO$. Namely,  
  \begin{equation}\label{BMOFunsNorms}
|f|_{BMO} \leq c, \quad \text{and}   \quad \,   |g|_{BMO} \leq c,  
\end{equation}
 with $c$  \textit{independent} of $\alpha$ (see Section \ref{Sec_BMOf}). Now we observe that if these functions have an angular dependance, for example  
   \begin{equation}\label{BMOFunsAng}
\til{f}(r,\theta)= \frac1{\beta(\alpha)}r^{\alpha} \sin(\theta), \quad \text{and}   \quad \,  \til{g}(r,\theta)= \frac{1}{\beta(\alpha)}r^{\alpha} \cos(\beta(\alpha)\log(r)) \sin(\theta). 
\end{equation}
Then, they become arbitrary  large in $BMO$. More precisely,  we have 
   \begin{equation}\label{BMOFunsAngNorms}
| \til{f} |_{BMO} \geq \frac1{\beta(\alpha)}  \rightarrow \infty  , \quad \text{and}   \quad \,   |\til{g}|_{BMO} \geq \frac1{\beta(\alpha)} \rightarrow \infty.     
\end{equation}
as $\alpha \rightarrow 0$, since $\lim_{\alpha \rightarrow 0} \beta(\alpha)=0$. See Section \ref{Sec_BMOf}) for more details.  

   \textbullet \textit{\, Ill-posedness mechanism:} 
We see from  \eqref{BMOFuns}, \eqref{BMOFunsNorms}, \eqref{BMOFunsAng}, and \eqref{BMOFunsAngNorms} that if a radial function such as \eqref{BMOFuns} or \eqref{BMOFunsNorms} develops an angular dependence, it will instantly become arbitrarily large in $BMO$. Hence, to obtain the ill-posedness result, we design an initial data for the 2d Euler equation in order for this to occur dynamically. We first observe that a compactly supported vorticity of the following form: 
$$
\omega_{odd}(r,\theta)=r^{\alpha}\sin(2\theta),
$$
 will generate a hyperbolic flow around the origin. Therefore, we choose   an initial data as:
$$
\omega(r,\theta)=h(r)+\omega_{odd}(r,\theta),
$$
where $h(r)$ is a radial function. Now if $h(r)$ is chosen carefully so that it does not destroy the hyperbolic dynamics from $\omega_{odd}(r,\theta)$, then the level sets of $h(r)$, which are initially circular, will develop an angular dependence and become elliptical. We choose $h(r)=g(r)$ from \eqref{BMOFuns} and we observe that initially
$$|\nabla u(\omega_{odd})|_{L^{\infty}} \big|_{t=0}\approx O(\frac{1}{\alpha}), \quad \text{and}   \quad  |\nabla u(h)|_{L^{\infty}} \big|_{t=0}\approx O(\ \frac{1}{\beta(\alpha)}).$$
Thus, by choosing $\beta(\alpha)$ significantly larger than $\alpha$, we can perturbatively treat the contribution from the vorticity $h(r)$ on a specific small time scale. Now, the reason we choose $h(r)$ to be $g(r)$ instead of, say, $f(r)$ in \eqref{BMOFuns} is slightly subtle. Although the gradient of the velocity field coming from $\omega_{odd}$, $\nabla u(\omega_{odd})$, is initially $O(\alpha^{-1})$, it gets immediately weakened by the nonlinear dynamics of the 2d Euler equations. Using data such as $f(r)$ from \eqref{BMOFuns} is not enough to obtain sufficient growth on our small time scale. Because of this, we use an oscillating radial function in order to make it more easier to leave $BMO$ when it gets advected by the hyperbolic flow from $\omega_{odd}$. Therefore, considering the described  initial data, it will  generate a vorticity of the following form
$$
\omega(t) \approx   \frac{1}{\beta(\alpha)}r^{\alpha} \cos(\beta(\alpha)\log(r)) \cos(2\theta),
$$   
which will instantly become arbitrarily large in $BMO$. This ill-posedness mechanism   is inspired by our  work \cite{SK}  on Besov spaces.  We observe that at the level of vorticity,  $BMO$ space and the Besov spaces $B^{0}_{\infty,q}$, for $1<q<\infty$, share  structural similarity. Namely, that if  the vorticity  is radial, then it can be rougher.

\textbullet \textit{\, Modeling the main dynamic through a leading order equation:} We use Elgindi's Biot-Savart law decomposition \cite{E} to simplify the 2d Euler dynamics. We derive a leading order model by considering  only the main terms of the velocity field identified in \cite{E}. (see Sections \ref{Sec_Prelim} and  \ref{Sec_BMO_LOM}). A similar strategy was used in the work with Elgindi \cite{ESK} to establish strong ill-posedness of bounded solutions to the Euler equation with Riesz forcing. We also used it recently to prove strong ill-posedness of 2d Euler in critical Besov spaces \cite{SK}. The ill-posedness of the leading order model will be done in Section \ref{Sec_LOM_ME}. Then, the result for full 2d Euler  will be done in Sections \ref{Sec_BMO_remainder_Est} and    \ref{Sec_Main}, using a  perturbative argument.  
 
\subsection{Notation}\label{SubSec_Notation}
We work primarily in the standard $(r,\theta)$ polar coordinates on $\mathbb{R}^2$. 
The $L^{\infty}$ and $C^{\alpha}$ H\"older norms will be written as $|\cdot|_{L^{\infty}}$ and $|\cdot|_{C^{\alpha}}$, respectively. In addition, we will use  H\"older-type norms due to Elgindi and Jeong  \cite{EJ3}, which are defined as follows:
 \begin{align}
|f|_{\mathring{C}^{1,\alpha}}=   |f|_{C^{\alpha}}+   |\dth  f|_{C^{\alpha}}+|r \dr f|_{C^{\alpha}}. 
 \end{align}
The Bounded Mean Oscillation ($BMO$) norm can be defined by: 
$$|f|_{BMO}=\sup_{B}\frac{1}{|B|}\Big|\int_{B} f-\frac{1}{|B|}\int_{B} f\Big|,$$
where $B$ ranges over balls, and $|B|$ denotes the measure of the ball. We will consider examples on $\mathbb{R}$ and on $\mathbb{R}^2$, so depending on the context, $B$  could be an interval or ball.

 \subsection{Organization}
 In Section \ref{Sec_Prelim} we state the  Biot-Savart law decomposition for the connivence  of the reader. Section \ref{Sec_BMOf} describes  the  H\"older approximation of the $BMO$ space. In Section \ref{Sec_BMO_LOM} we derive what we call the leading order model for 2d Euler equation.  In Section \ref{Sec_LOM_ME}, we prove the main estimates for the leading order model.  Section \ref{Sec_LOM_UE}  provides useful estimates for the model which we will use  in Section \ref{Sec_BMO_remainder_Est} to estimate the remainder. Finally, in Section \ref{Sec_Main},  we prove our main result.

 \section{Preliminary Results:}\label{Sec_Prelim}
In this section, we describe the Biot-Savart law decomposition by  Elgindi \cite{E}. This will be used in Section \ref{Sec_BMO_LOM} to derive our leading order model.

\begin{proposition}\label{ElgindiDecomp}(Elgindi \cite{E}: Biot-Savart law decomposition) Consider compactly supported vorticity $\omega \in C^{\alpha}(\mathbb{R}^{2})\cap {\mathring{C}^{1,\alpha}}(\mathbb{R}^{2})$. Writing  the stream function as: $$\psi(r,\theta)=r^2\Psi(r,\theta)$$ 
Then we have the following  decomposition: 
\begin{align*}
 \Psi(\omega)(r,\theta)=L^{s}_{12}(\omega)(r) \sin(2\theta)+L^{c}_{12}(\omega)(r) \cos(2\theta)+\Psi_r(\omega)(r,\theta)
 \end{align*}
 where
 \begin{align*}
 L^s_{12}(\omega)(r)&= \frac{1}{4\pi} \int_{r}^{\infty} \int_0^{2\pi} \frac{\omega(s,\theta)}{s} \sin(2\theta) \, d\theta  ds,\\
  L^c_{12}(\omega)(r)&=\frac{1}{4\pi} \int_{r}^{\infty} \int_0^{2\pi} \frac{\omega(s,\theta)}{s}  \cos(2\theta) \, d\theta ds,
 \end{align*}
  and with  the remainder term $\Psi_r(\omega)(r,\theta)$ satisfying 
  $$|\Psi_r(\omega)|_{C^{\alpha}}\leq c |\omega|_{C^{\alpha}\cap {\mathring{C}^{1,\alpha}}}, \, \text{and} \,\,\, |r\dr \Psi_r(\omega)|_{C^{\alpha}} \leq c |\omega|_{C^{\alpha}\cap {\mathring{C}^{1,\alpha}}}, $$
for a constant  $c$ is independent of $\alpha$. 
\end{proposition}
\begin{remark}
The Biot-Savart law decomposition in \cite{E} was done for the 3d Euler equation. It is  simpler in two dimensions and it has different versions, including an $L^{\infty}$ version, see \cite{Eremarks} and Lemma 5.1 in \cite{EJ2}, and a weighted $L^{2}$ energy version, see Theorems 4.23 and 4.24 in  \cite{DE}. The above decomposition we use is obtained by combining these mentioned results with Theorem 2 in \cite{EJ3}. We mention that one may view the decomposition as a generalization  of  Kiselev and \v{S}ver\'{a}k Key Lemma in \cite{KS}. 
\end{remark}

 \begin{remark}\label{Remark_psi_remainder}
In Section \ref{Sec_BMO_remainder_Est}, when obtaining the remainder estimates, we will use the lower order part of the stream function $\psi_r$ defined as:
\begin{align*}
\psi_r :&= \psi(r,\theta) - r^2L^{s}_{12}(\omega)(r) \sin(2\theta) - r^2L^{c}_{12}(\omega)(r) \cos(2\theta). 
\end{align*}
This term gives $C^{\alpha}$ H\"older bounds with constant independent of $\alpha$.
\end{remark}

\section{$BMO$ space using H\"older-type Functions }\label{Sec_BMOf}
In this section, we approximate functions in the $BMO$ space by H\"older-type  functions. This is done by coupling the H\"older regularity exponent with an appropriate rescaling, which gives quantitative estimates for the $BMO$ norm in terms of the regularity exponent. This now allows us to study unbounded functions using a rescaling of H\"older-type functions for which the 2d Euler equation is well-posed. Further, using this approach, we identify two families of functions in $BMO$, observing  that a function in  $BMO$ can be unbounded with mild oscillation, see Lemma \ref{TypeII}. This observation will be leveraged in the ill-posedness mechanism, see Section \ref{Sec_LOM_ME}. We remark that this approach is inspired by a similar type of approximation scheme  introduced in our previous work \cite{SK} on the strong ill-posedness of the 2d Euler equation in Besov spaces.  There, it was used to approximate bounded vorticity, whereas here we extend it to the  $BMO$ space to now approximate unbounded vorticity.

For this section, the results hold as long as the rescaling $\beta(\alpha)$ factor satisfies  the conditions $\lim_{\alpha \rightarrow 0} \beta(\alpha)=0$, and $\frac{\alpha}{\beta(\alpha)}=O(1)$.  For example, $\beta(\alpha)=\alpha$ would work. However, for the purpose of the ill-posedness of the 2d Euler equation, we will need to choose a specific scaling which is going to be $\beta(\alpha)= |\log(|\log(\alpha)|^{1+\eta})|^{-1}$ for some $\eta>0$. 

This section is organized as follows: In Lemma \ref{TypeI}, we give one  type of functions in $BMO$. Then,  in Lemma \ref{TypeII}, we give the second type of oscillating functions in $BMO$. Finally, in Lemma \ref{BMOexLow}, we give an example of how these function can become arbitrarily large in $BMO$.

\begin{lemma}\label{TypeI} Let  $0<\alpha<1$ and consider the   function $f:[0,1] \rightarrow \mathbb{R} $ defined as follows: 
$$f(r)=   
 \frac1{\beta(\alpha)}r^{\alpha},
 $$
such that  $\lim_{\alpha \rightarrow 0} \beta(\alpha)=0$, and $\frac{\alpha}{\beta(\alpha)}=O(1)$. Then we have the following:  
$$
 |f|_{BMO} \leq  c,
$$
 where $c$ is independent of $\alpha$. 
\end{lemma}

\begin{proof} Here we follow the strategy of proving that a function belongs to $BMO$  that can be found in Grafakos \cite{Gr} Section 3, and we adjust it to our case.  Given an interval $(-h+r_0,r_0+h)$ for  $h>0$, it suffices to show that there exist a constant $c_{r_0,h}$ such that 
\begin{align}\label{IntMain}
\frac{1}{2h}\int_{r_0-h}^{r_0+h}| \frac1{\beta(\alpha)}r^{\alpha}-c_{r_0,h}|dr \leq c
\end{align}
where $c$ is independent of $r_0,h$, and $\alpha$. We will consider two cases when $r_0>2h$, and when $r_0\leq 2h$. First, we consider the case when $r_0>2h$. Here we take $c_{r_0,h}= \frac1{\beta(\alpha)}r_0^{\alpha}$.  Thus, we have
\begin{align}\label{Int}
\frac{1}{2h}\int_{r_0-h}^{r_0+h}| \frac1{\beta(\alpha)}r^{\alpha}-\frac1{\beta(\alpha)}r_0^{\alpha}|dr=\frac{1}{2h}\int_{r_0}^{r_0+h} ( \frac1{\beta(\alpha)}r^{\alpha}-\frac1{\beta(\alpha)}r_0^{\alpha} )dr + \frac{1}{2h}\int_{r_0-h}^{r_0} (-  \frac1{\beta(\alpha)}r^{\alpha}+\frac1{\beta(\alpha)}r_0^{\alpha})dr. 
\end{align}
For the first integral, we obtain 
\begin{align}\label{Int1}
\frac{1}{2h}\int_{r_0}^{r_0+h} ( \frac1{\beta(\alpha)}r^{\alpha}-\frac1{\beta(\alpha)}r_0^{\alpha} )dr =\frac{1}{2h}( \frac1{\beta(\alpha)}\frac{1}{\alpha+1}(r_0+h)^{\alpha+1}- \frac1{\beta(\alpha)}\frac{1}{\alpha+1}r_0^{\alpha+1})   - \frac{1}{2h}  \frac1{\beta(\alpha)}r_0^{\alpha} h.  
\end{align}
For the second integral, we have 
\begin{align}\label{Int2}
 \frac{1}{2h}\int_{r_0-h}^{r_0} ( -\frac1{\beta(\alpha)}r^{\alpha}+\frac1{\beta(\alpha)}r_0^{\alpha})dr = \frac{1}{2h} ( -\frac1{\beta(\alpha)}\frac{1}{\alpha+1}r_0^{\alpha+1}+\frac1{\beta(\alpha)}\frac{1}{\alpha+1}(r_0-h)^{\alpha+1})+ \frac{1}{2h}   \frac1{\beta(\alpha)}r_0^{\alpha}h.
\end{align}
 From \eqref{Int},\eqref{Int1}, and  \eqref{Int2}, we observe that 
\begin{align*} 
 \frac{1}{2h}\int_{r_0-h}^{r_0+h}| \frac1{\beta(\alpha)}r^{\alpha}-\frac1{\beta(\alpha)}r_0^{\alpha}|dr&=\frac{1}{2h} \frac1{\beta(\alpha)}\frac{1}{\alpha+1} \Big( (r_0+h)^{\alpha+1}-2 r_0^{\alpha+1} + (r_0-h)^{\alpha+1} \Big)  \\
 &\leq  \frac{1}{2h} \frac1{\beta(\alpha)}\frac{\alpha }{\alpha+1}cr_0^{\alpha+1} \frac{h^2}{r_0^{2}}=   \frac{c}{\beta(\alpha)}\frac{\alpha }{\alpha+1}r_0^{\alpha-1}h,
 \end{align*}
where, for the inequality, we used that $r_0>2h$ and the properties of the second finite difference. Using again that    that $r_0>2h$ with the  assumption that $\frac{\alpha}{\beta(\alpha)}\leq c$, we obtain  
$$
 \frac{1}{2h}\int_{r_0-h}^{r_0+h}| \frac1{\beta(\alpha)}r^{\alpha}-\frac1{\beta(\alpha)}r_0^{\alpha}|dr \leq  \frac{\alpha}{\beta(\alpha)}\frac{h}{r_0}r_0^{\alpha} \leq c,
$$
where $c$ is independent of  $r_0,h$, and $\alpha$. Hence, we have \eqref{IntMain} for the first case when $r_0>2h$. Now we consider the second case, when $r_0\leq 2h$. In this case, it suffices to show that there exist a constant $c_h$ such that  
$$
 \frac{1}{h}\int_{0}^{h}| \frac1{\beta(\alpha)}r^{\alpha}-c_{h}|dr \leq c,
 $$
for some $c$ independent  of $h$ and $ \alpha$. We set $c_{h}=\frac1{\beta(\alpha)}h^{\alpha}$ and observe
$$
 \frac{1}{h}\int_{0}^{h}| \frac1{\beta(\alpha)}r^{\alpha}-c_{h}|dr = \frac{1}{h}\int_{0}^{h}| \frac1{\beta(\alpha)}r^{\alpha}-\frac1{\beta(\alpha)}h^{\alpha}|dr=\frac{1}{h}(\frac1{\beta(\alpha)}h^{\alpha+1}- \frac1{\beta(\alpha)}\frac{1}{\alpha+1}h^{\alpha+1}).
 $$
Therefore,  using the assumption  that $\frac{\alpha}{\beta(\alpha)}\leq c$, we have 

$$
 \frac{1}{h}\int_{0}^{h}| \frac1{\beta(\alpha)}r^{\alpha}-c_{h}|dr =\frac{h^{\alpha+1}}{h}\frac{1}{\beta(\alpha)}(1-\frac{1}{1+\alpha}) \leq\frac{h^{\alpha+1}}{h}\frac{c\alpha}{\beta(\alpha)} \leq ch^{\alpha} \leq c,$$
where $c$ is   independent  of $h$ and $\alpha$.  This completes the proof.

\end{proof}
\begin{lemma}\label{TypeII}  Let  $0<\alpha<1$ and consider a   function $g:[0,\frac{1}{2}] \rightarrow \mathbb{R} $ defined as follows: 
$$g(r)=   
 \frac{1}{\beta(\alpha)}r^{\alpha} \cos(\beta(\alpha)\log(r)), \,\, \text{or}\,\,\,     g(r)=   
 \frac{1}{\beta(\alpha)}r^{\alpha} \sin(\beta(\alpha)\log(r)),
 $$
such that  $\lim_{\alpha \rightarrow 0} \beta(\alpha)=0$, and $\frac{\alpha}{\beta(\alpha)}=O(1)$. Then we have the following:  
$$
 |g|_{BMO} \leq  c
$$
 where $c$ is independent of $\alpha$. 

\end{lemma}
\begin{remark}We note that  mild oscillation (the  $\beta(\alpha)\log(r)$ factor) is needed. For instance, the function  $\frac{1}{\beta(\alpha)}r^{\alpha} \cos(\log(r))$ is not in $BMO$. 
\end{remark} 

\begin{proof}
The argument here is similar to the previous lemma, so we will only consider the following case, where we will show that there exist a constant $c_{h}$ such that
$$
\frac{1}{h}\int_{0}^{h}| \frac1{\beta(\alpha)}r^{\alpha} \cos(\beta(\alpha)\log(r))-c_{h}|dr \leq c,
$$
for some constant $c$ independent  of $h$ and $\alpha$. We set $c_{h}=\frac{1}{\beta(\alpha)}h^{\alpha} \cos(\beta(\alpha)\log(h))$ and obtain  
\begin{align} 
\frac{1}{h}\int_{0}^{h}| \frac1{\beta(\alpha)}r^{\alpha} \cos(\beta(\alpha)\log(r))&-\frac{1}{\beta(\alpha)}h^{\alpha} \cos(\beta(\alpha)\log(h))|\\
&\leq    \frac{1}{h}\int_{0}^{h}| \frac1{\beta(\alpha)}r^{\alpha} \cos(\beta(\alpha)\log(r))-\frac{1}{\beta(\alpha)}h^{\alpha} \cos(\beta(\alpha)\log(r))|dr \label{IntLem2A} \\
&+ \frac{1}{h}\int_{0}^{h}| \frac1{\beta(\alpha)}h^{\alpha} \cos(\beta(\alpha)\log(r))-\frac{1}{\beta(\alpha)}h^{\alpha} \cos(\beta(\alpha)\log(h))|dr\label{IntLem2B} 
\end{align} 
\noindent  For the integral in \eqref{IntLem2A}  we have 
\begin{align*}
      \frac{1}{h}\int_{0}^{h}| \frac1{\beta(\alpha)}r^{\alpha} \cos(\beta(\alpha)\log(r))-\frac{1}{\beta(\alpha)}h^{\alpha} & \cos(\beta(\alpha)\log(r))|dr\\ &\leq   \frac{1}{h}\int_{0}^{h}| \frac1{\beta(\alpha)}r^{\alpha}-\frac{1}{\beta(\alpha)}h^{\alpha} |  |\cos(\beta(\alpha)\log(r))| dr \\
       &\leq   \frac{1}{h}  \frac{1}{\beta(\alpha)}  \big( h^{\alpha+1}-  \frac{1}{1+\alpha} h^{\alpha+1}  \big) =   \frac{h^{\alpha}}{\beta(\alpha)}  \big( 1-  \frac{1}{1+\alpha}  \big) \\
       & \leq     c \frac{\alpha}{\beta(\alpha)}   h^{\alpha}.  
      \end{align*} 
Using the assumption that $\frac{\alpha}{\beta(\alpha)}\leq c$, we have 
  \begin{align}\label{IntLem2AEst}
      \frac{1}{h}\int_{0}^{h}| \frac1{\beta(\alpha)}r^{\alpha} \cos(\beta(\alpha)\log(r))-\frac{1}{\beta(\alpha)}h^{\alpha} \cos(\beta(\alpha)\log(r))|dr &\leq c,
      \end{align} 
where $c$ is independent of $\alpha$, and $h$. Next, we consider  the \eqref{IntLem2B} integral 
       \begin{align}\label{IntLem2BEst1}
 \frac{1}{h}\int_{0}^{h}| \frac1{\beta(\alpha)}h^{\alpha} \cos(\beta(\alpha)\log(r))-\frac{1}{\beta(\alpha)}h^{\alpha} &\cos(\beta(\alpha)\log(h))|dr\\ &= \frac{1}{h}\int_{0}^{h} \frac{h^{\alpha}}{\beta(\alpha)}  |   \cos(\beta(\alpha)\log(r))-\cos(\beta(\alpha)\log(h))|dr    
 \end{align} 
Observe that we have 
\begin{align}\label{IntLem2BEst2}
|\cos(\beta(\alpha)\log(r))-\cos(\beta(\alpha)\log(h))|\leq \int_r^h \frac{\beta(\alpha)}{s} ds \leq \beta(\alpha) |\log(\frac{r}{h})|
 \end{align} 
Thus, using \eqref{IntLem2BEst1} and \eqref{IntLem2BEst2} we obtain 
$$
 \frac{1}{h}\int_{0}^{h}| \frac1{\beta(\alpha)}h^{\alpha} \cos(\beta(\alpha)\log(r))-\frac{1}{\beta(\alpha)}h^{\alpha} \cos(\beta(\alpha)\log(h))|dr  \leq c h^{\alpha} \int_0^{1}|\log(r)| dr \leq c h^{\alpha} $$
 Therefore, this gives
  \begin{align}\label{IntLem2BEst}
 \frac{1}{h}\int_{0}^{h}| \frac1{\beta(\alpha)}h^{\alpha} \cos(\beta(\alpha)\log(r))-\frac{1}{\beta(\alpha)}h^{\alpha} \cos(\beta(\alpha)\log(h))|dr \leq    c     
 \end{align} 
where $c$ is independent of $\alpha$. Thus, from  \eqref{IntLem2AEst} and \eqref{IntLem2BEst} combined with \eqref{IntLem2A} and  \eqref{IntLem2B} we obtain that
$$
\frac{1}{h}\int_{0}^{h}| \frac1{\beta(\alpha)}r^{\alpha} \cos(\beta(\alpha)\log(r))-\frac{1}{\beta(\alpha)}h^{\alpha} \cos(\beta(\alpha)\log(h))|\leq c
$$
with $c$ independent  of $h$ and $\alpha$. Thus, completing the proof.  

\end{proof}


\begin{lemma}\label{BMOexLow}  Let  $0<\alpha<1$ and consider the following vortices on $\mathbb{R}^2$: 
$$\omega(r,\theta)= f(r)a(\theta) \quad \text{or} \quad  \omega(r,\theta)=g(r)a(\theta), 
 $$
where $f(r), g(r)$ are defined as in Lemmas \ref{TypeI} and \ref{TypeII}, respectively, and $a(\theta)$ is any smooth mean zero angular function. Then we have the following:  
$$
 |\omega|_{BMO} \geq \frac{c}{\beta(\alpha)},
$$
where $c$ is independent of $\alpha$. Since $\beta(\alpha) \rightarrow 0 $ as $\alpha \rightarrow 0$, we have $ |\omega|_{BMO} \rightarrow \infty$.

\end{lemma}

\begin{proof}
We prove the result for the case $ \omega(r,\theta)=g(r)a(\theta)$, since the other case is similar. 
$$
 |\omega|_{BMO} \geq \frac{1}{\pi h^2}\int_0^{h} \int_0^{2\pi} \bigg|   \frac1{\beta(\alpha)}r^{\alpha} \cos(\beta(\alpha)\log(r)) a(\theta) -\frac{1}{\pi h^2}\int_{0}^{h} \int_0^{2\pi}  \frac1{\beta(\alpha)}s^{\alpha} \cos(\beta(\alpha)\log(s)) h(\sigma) s \,d\sigma  ds    \bigg| r  d\theta \, dr. 
$$
Since  $h$ is  mean zero angular function,  the inner integral is zero. Thus, we obtain 

$$
 |\omega|_{BMO} \geq \frac{1}{\pi h^2}\int_0^{h} \int_0^{2\pi} \bigg|   \frac1{\beta(\alpha)}r^{\alpha} \cos(\beta(\alpha)\log(r)) a(\theta)    \bigg| r  d\theta  dr 
$$
Using again that $a(\theta)$ is mean  zero, we can find an interval $(\theta_1,\theta_2)$ for which $a(\theta)>c>0$. Hence,   
\begin{align}
 |\omega|_{BMO} &\geq \frac{1}{\pi h^2}\int_0^{h} \int_{\theta_1}^{\theta_2} \bigg|   \frac1{\beta(\alpha)}r^{\alpha} \cos(\beta(\alpha)\log(r)) a(\theta)    \bigg| r  d\theta \, dr \label{BMOest1}\\
      & \geq\frac{c}{\pi h^2}\int_0^{h}   \frac1{\beta(\alpha)}r^{1+\alpha} |\cos(\beta(\alpha)\log(r))|   \label{BMOest2}     dr 
 \end{align}
Changing variable to $u=\beta(\alpha) \log(r)$, we have 
\begin{align}\label{BMOestInt}
 \frac{c}{\pi h^2}\int_0^{h}   \frac1{\beta(\alpha)}r^{\alpha+1} |\cos(\beta(\alpha)\log(r))|         dr=  \frac{c}{\pi h^2} \frac1{\beta^2(\alpha)}  \int_{\beta(\alpha) |\log(h)|}^{\infty}  e^{-\frac{u}{\beta(\alpha)}(\alpha+2)}|\cos(u)|       du  
 \end{align}
We observe that for any positive integer $k$, we obtain 

\begin{align*}
  \int_{\frac{3}{4}\pi+k\pi}^{\infty}  e^{-\frac{u}{\beta(\alpha) } (\alpha+2)}  | \cos( u)  |         du   &\geq \sum_{j=k}^{\infty}  \int_{\frac{3}{4}\pi+j\pi}^{\frac{5}{4}\pi+j\pi}  \frac{1}{\sqrt{2}}e^{-\frac{u}{\beta(\alpha) } (\alpha+2)}          du  \\
      &\geq \sum_{j=k}^{\infty}     c e^{-\frac{\frac{5}{4}\pi+j\pi}{\beta(\alpha) } (\alpha+2)}        =  c \frac{e^{-\frac{\frac{5}{4}\pi+k\pi}{\beta(\alpha) } (\alpha+2)}}{1-e^{-\frac{\pi}{\beta(\alpha) } (\alpha+2)}  }
      \end{align*} 
 Now since  $\beta(\alpha)\rightarrow 0$ as $\alpha \rightarrow 0$, this gives 
  
 $$   \int_{\frac{3}{4}\pi+k\pi}^{\infty}  e^{-\frac{u}{\beta(\alpha) } (\alpha+2)}  | \cos( u)  |         du     \geq c\beta(\alpha) e^{-\frac{\frac{5}{4}\pi+k\pi}{\beta(\alpha) } (\alpha+2)}.$$
  Thus, by choosing $\alpha$ appropriately, we have 
\begin{align}\label{BMOestInt2} 
 \int_{\beta(\alpha) |\log(h)|}^{\infty}    e^{-\frac{u}{\beta(\alpha) } (\alpha+2)}  | \cos( u)  |         du     \geq  c\beta(\alpha)  e^{-|\log(h)| (\alpha+2)}= c\beta(\alpha) h^2.    
 \end{align}
 Therefore, from \eqref{BMOest2}, \eqref{BMOestInt}, and \eqref{BMOestInt2}, we have our desired result
 $$
 |\omega|_{BMO}  \geq    \frac{c}{\beta(\alpha)}.
$$
\end{proof}
\section{Leading Order Model}\label{Sec_BMO_LOM}
We derive the leading order model for the 2d Euler equation. 
We remark that in our work with Elgindi \cite{ESK} on the Euler equation  with Riesz forcing, a similar model was derived in order to prove strong ill-posedness in $L^{\infty}$. The model we consider here is the same model we used recently in \cite{SK} to prove strong ill-posedness of the 2d Euler in critical Besov spaces.  

The vorticity form of the 2d Euler equation  in polar coordinates is:  
\begin{equation}
\partial_t \omega + \frac{1}{r} \, \partial_r \psi \, \partial_\theta \omega 
- \frac{1}{r} \, \partial_\theta \psi \, \partial_r \omega = 0,
\end{equation}
where $\psi$ is the stream function which solves $\Delta \psi = \omega$. 
Using the Biot Savart law decomposition from Proposition \ref{ElgindiDecomp}, by setting $\psi(r,\theta)=r^2\Psi(r,\theta)$, the stream function can be written as:
 $$
 \Psi(\omega)(r,\theta)=L^{s}_{12}(\omega)(r) \sin(2\theta)+L^{c}_{12}(\omega)(r) \cos(2\theta)+\Psi_r(\omega)(r,\theta),
 $$
 where we recall that the operators $L^s_{12}$ and $L^c_{12}$ are:
\begin{align}\label{BMOL12op}
 L^{s}_{12} (\omega)(r) &= \frac{1}{4\pi} \int_r^{\infty} \int_0^{2\pi} \frac{\omega(s,\theta)}{s}  \sin(2\theta) d\theta   ds, \, \text{and} \,\, 
 L^{c}_{12} (\omega)(r) = \frac{1}{4\pi} \int_r^{\infty} \int_0^{2\pi} \frac{\omega(s,\theta)}{s}   \cos(2\theta) d\theta ds. 
\end{align}
Then, the model we consider is:
\begin{equation} \label{LOMeq}
 \dt \Omega+  2 L^s_{12}(\Omega) \sin(2\theta)  \dth \Omega- 2 rL^s_{12}(\Omega) \cos(2\theta) \dr \Omega =0,
\end{equation} 
with initial data 
\begin{equation} \label{LOMeqOmega_0}
\Omega_0(r,\theta) =\frac{1}{\beta(\alpha)}r^{\alpha} \cos(\beta(\alpha)\log(r))+  r ^{\alpha} \sin(2\theta),
\end{equation} 
compactly supported on the ball of radius 1. From now on  $\beta(\alpha)=\frac{1}{\log(|\log(\alpha)|^{1+\eta})}$ for some $\eta>0$.

\section{Leading Order Model Main Estimates}\label{Sec_LOM_ME}
In this section, the main estimates for the leading order model are proved. Proposition \ref{L12_Up_Low} provides an upper and lower bound on the main term of the gradient of the velocity field. These will be used in  Proposition \ref{LOM_BMO} to prove the strong ill-posedness of the model in $BMO$.  

  \begin{proposition} \label{L12_Up_Low} 
   Let $\Omega(t)$ be a solution to the leading order model \eqref{LOMeq}   with initial data \eqref{LOMeqOmega_0}. Then,   there exist constants $c$ independent of $\alpha$ such that we have the follow estimates
 
\noindent  The upper bound estimate: 
\begin{equation} \label{L12_Up_est}
 |L^s_{12}(\Omega)(t)|_{L^{\infty}} \leq \frac{c}{\alpha} ,
\end{equation} 
and the  lower   bound estimate:

\begin{equation} \label{L12_Low_est}
\sup_{0\leq r \leq 1}| \int_0^{t}L^s_{12}(\Omega) d\tau|\geq  c \log(1+c\frac{t}{\alpha}),
\end{equation} 
\noindent on the time scale  $0\leq t \leq T(\alpha)=c \alpha \log(c |\log(\alpha)|) $, with c independent of $\alpha$. 
   
   \end{proposition}
   
   \begin{remark}
In \cite{ESK} with Elgindi, similar type of estimates was derived, and also recently in \cite{SK}. 
\end{remark}

   \begin{proof} We will first prove the upper bound  \eqref{L12_Up_est} and then we will use it to prove the lower bound  \eqref{L12_Low_est}. We  study the flow map of the leading order model, which will be used for both estimates. Recall from \eqref{LOMeq} that  the  leading order model is:  
$$
 \dt \Omega+  2 L^s_{12}(\Omega) \sin(2\theta)  \dth \Omega- 2 rL^s_{12}(\Omega) \cos(2\theta) \dr \Omega =0.
 $$
We consider the characteristics of the model and  observe that we can write a coupled system of ODEs for the angular and radial trajectories of the flow map. To shorten the notation, we will write  for now  
\begin{equation*} \Phi(t,r,\theta)=(\Phi_r (t,r,\theta),\Phi_{\theta}(r,\theta)))=(r(t),\theta(t))
\end{equation*} 
 where  we have   $\Phi( 0,r,\theta) =(r(0),\theta(0))=(r,\theta)$. Therefore, we can obtain  the following:  
 \begin{align}
\dot{r}(t)&=  - 2r(t) L^s_{12}(\Omega) \cos(2\theta(t))\label{RadialT_Up_Low}, \\
\dot{\theta}(t)&=2 L^s_{12}(\Omega)  \sin(2\theta(t))\label{AnguarT_Up_Low}.
\end{align}
Hence, this gives 
 $$
\frac{dr}{d\theta }=-\frac{  r   \cos(2\theta)}{     \sin(2\theta)}.
$$
which we can solve and obtain  
\begin{equation}\label{FlowMap_Up_Low}
 (r(t))^2 |\sin(2\theta(t))|   =     r^2  |\sin(2\theta)| .
\end{equation} 
 We observe that given  a point  starting at  $(r,\theta)$, equation \eqref{FlowMap_Up_Low} corresponds to the trajectories of a hyperbolic flow around the origin, which is also satisfied by the  inverse flow map with the trajectories moving backward. With a slight abuse of notation, from now on,  we will write 
 $$\Phi^{-1}(t,r,\theta)=(\Phi^{-1}_r (t,r,\theta),\Phi^{-1}_{\theta}(r,\theta)))=(r(t),\theta(t)),$$where   $\Phi^{-1}( 0,r,\theta) =(r(0),\theta(0))=(r,\theta)$. In addition, to shorten the notation, we will write
\begin{equation}\label{XFlow_Up_Low}
  X(t,r,\theta):=2\int_0^t L_{12}(r(\tau))\cos(2\theta(\tau)) d\tau. 
\end{equation} 
 Therefore, from \eqref{RadialT_Up_Low}, we  have 
 
  $$r(t)=re^{X(t,r,\theta)}.$$ 
 Recall from $\eqref{LOMeqOmega_0}$ that the initial data  is:  
 $$\Omega_0(r,\theta) =\frac{1}{\beta(\alpha)}r^{\alpha} \cos(\beta(\alpha)\log(r))+  r ^{\alpha} \sin(2\theta).$$
Since the vorticity is being transported, we can write the solution as 
 \begin{align}
 \Omega(t,r,\theta)= \Omega_0(\Phi^{-1}(t,r,\theta))&= \frac{1}{\beta(\alpha)} (r(t))^{\alpha}  \cos(\beta(\alpha)\log(r(t)))+ (r (t) )^{\alpha} \sin(2\theta(t)), \\
 &= \frac{1}{\beta(\alpha)} r^{\alpha}  e^{\alpha X(t)} \cos(\beta(\alpha)\log(r e^{X(t)})+ r^{\alpha}  e^{\alpha X(t)} \sin(2\theta(t)) \label{OmgeSol_Up_Low}
\end{align}
 Now we can  begin proving our desired bounds, so we start with the upper bound. 
 
\noindent \textbf{Upper Bound \eqref{L12_Up_est}:} Here we use a bootstrap argument. More precisely, we assume that there is a constant $c$, independent of $\alpha$, such that 
\begin{equation}\label{Bootstrap_H_Up_es}
  |L^s_{12}(\Omega)(t)|_{L^{\infty}} \leq   \frac{10c}{\alpha},  
\end{equation} 
Then, we show that
  \begin{equation}\label{Bootstrap_C_Up_es}
  |L^s_{12}(\Omega)(t)|_{L^{\infty}} \leq    \frac{c}{\alpha}.  
\end{equation} 
  Recall from \eqref{BMOL12op} that  the operator  $ L^{s}_{12}$ is defined as follows: 
 
  $$
 L^{s}_{12} (\Omega)(r)=c \int_r^{\infty} \int_0^{2\pi}  \frac{ 1}{s} \Omega(t,s,\theta) \sin(2\theta) d \theta \, ds.
$$
For the angular integral, due to symmetry, it suffices to consider the  integral on $[0,\frac{\pi}{2}]$. Further, for radial integral, since the initial data has compact support, we can consider the integral on $[0,1]$. Thus, we have 
$$
 L^{s}_{12} (\Omega)(r)=c \int_r^{1} \int_0^{\frac{\pi}{2}}  \frac{ 1}{s} \Omega(t,s,\theta) \sin(2\theta) d \theta \, ds
$$
From \eqref{OmgeSol_Up_Low}, we can write $L_{12}$ as follows: 
\begin{equation}\label{L12_Up_Low_OP}
  L_{12}^{s}(\Omega)=c\int_{r}^{1}\int_{0}^{\frac{\pi}{2}}  \frac{1}{s}   \Big( \frac{1}{\beta(\alpha)}  s^{\alpha}    e^{\alpha X(t)} \cos(\beta(\alpha)\log(s e^{X(t)}))+ s^{\alpha}  e^{\alpha X(t)} \sin(2\theta(t)) \Big)    \sin(2\theta) d\theta ds. 
\end{equation}   For the first term, we observe that since $X(t)=2\int_0^t(L_{12}(s(\tau))\cos(2\theta(\tau)))d\tau$ is  even in $\theta$, we have 
\begin{equation}\label{Term1_Up_Low} 
\int_{r}^{1}\int_{0}^{\frac{\pi}{2}}  \frac{1}{s}  \frac{1}{\beta(\alpha)}  s^{\alpha}    e^{\alpha X(t)} \cos(\beta(\alpha)\log(s e^{X(t)}))   \sin(2\theta) d\theta ds=0
\end{equation} 
Hence, $L_{12}$ reduces to
\begin{equation}\label{L12_Up_Low_2}
  L_{12}^{s}(\Omega)=c\int_{r}^{1}\int_{0}^{\frac{\pi}{2}}  \frac{  s^{\alpha}  }{s}  e^{\alpha X(t)} \sin(2\theta(t))      \sin(2\theta) d\theta ds. 
\end{equation} 
 Now we use the bootstrap assumption. Namely,  from \eqref{Bootstrap_H_Up_es} and  \eqref{XFlow_Up_Low}, we have 
$$
|X(t)|_{L^{\infty}} \leq  \frac{10c}{\alpha} t \implies \alpha |X(t)|_{L^{\infty}} <   10c  t <1,
$$
on the time interval $0\leq t \leq T(\alpha)=c \alpha \log(c|\log(\alpha)|) $. Thus, the term $ |e^{\alpha X(t)}| $ can be made uniformly bounded on the time scale $ T(\alpha)$. Therefore, from \eqref{L12_Up_Low_2}, we have 
\begin{equation} 
  |L_{12}^{s}(\Omega)| \leq  c \int_{r}^{1}\int_{0}^{\frac{\pi}{2}}\frac{ s^{\alpha} }{s}    | \sin(2\theta(t))|    |\sin(2\theta)| d\theta ds \leq  c \int_{r}^{1}\int_{0}^{\frac{\pi}{2}}\frac{    s^{\alpha}}{s}      d\theta ds \leq \frac{c}{\alpha},
\end{equation}
and by adjusting the constant $c$, if necessary,  we have the desired upper bound:  
\begin{equation} 
 | L_{12}^{s}(\Omega)|_{L^{\infty}} \leq  \frac{c}{\alpha} , 
\end{equation} 
which completes the bootstrap argument \eqref{Bootstrap_C_Up_es}, so this gives \eqref{L12_Up_est}.

\noindent \textbf{Lower Bound:} To obtain the lower bound, from \eqref{L12_Up_Low_2}, we have 

\begin{equation}\label{L12_Up_Low_3}
  L_{12}^{s}(\Omega)=c\int_{r}^{1}\int_{0}^{\frac{\pi}{2}}  \frac{  s^{\alpha}  }{s}  e^{\alpha X(t)} \sin(2\theta(t))      \sin(2\theta) d\theta ds. 
\end{equation}
From \eqref{FlowMap_Up_Low},  $ L_{12}^{s}(\Omega)$ can be written as: 
 \begin{equation}\label{L12_Up_Low_4}
     L_{12}^{s}(\Omega)(t)= c \int_{r}^{1}\int_{0}^{\frac{\pi}{2}}\frac{    s^{\alpha} }{s}   e^{(\alpha-1) X(t)}            \sin^2(2\theta) d\theta ds. 
 \end{equation} 
We  take the time  derivative 
 \begin{align}\label{L12_Up_Low_5}
     \partial_t     L_{12}^{s}(\Omega)(t)   &= c \int_{r}^{1}\int_{0}^{\frac{\pi}{2}}\frac{    s^{\alpha} }{s}   e^{(\alpha-1) X(t)}      2(\alpha-1) X'(t)        \sin^2(2\theta) d\theta ds\\ 
       &=c \int_{r}^{1}\int_{0}^{\frac{\pi}{2}} \frac{     s^{\alpha} }{s}   e^{(\alpha-1) X(t)}   2(\alpha-1)         L^s_{12}(\Omega(t)) \cos(2\theta(t) )      \sin^2(2\theta) d\theta ds
 \end{align}
Hence, we have the following

 \begin{equation}\label{dtL12_Up_Low}
    \partial_t     L_{12}^{s}(\Omega)(t)\geq - c \int_{r}^{1}\int_{0}^{\frac{\pi}{2}}\frac{  s^{\alpha} }{s}   e^{(\alpha-1) X(t)}      2    L^s_{12}(\Omega(t))        \sin^2(2\theta) d\theta ds. 
 \end{equation} 
In addition, by taking the radial derivative of \eqref{L12_Up_Low_4} we obtain  
 \begin{equation}\label{drL12_Up_Low}
    \dr L_{12}^{s}(\Omega)(t)= - c\int_{0}^{2\pi}\frac{  r^{\alpha} }{r} e^{(\alpha-1) X(t)}     \sin^2(2\theta) d\theta ds. 
 \end{equation}
  Thus, combining \eqref{dtL12_Up_Low} and \eqref{drL12_Up_Low} we arrive at the  Riccati-type equation for $L_{12}^{s}$ 
    \begin{equation}\label{L12Riccat_Up_Low}
   \partial_t  L_{12}^{s}(\Omega)(t)\geq  \int_{r}^{1} 2    \partial_{s} L_{12}^{s}(\Omega)  L^s_{12}(\Omega(\tau)) = \int_{r}^{1}    \partial_{s} \Big ( L_{12}^{s}(\Omega)  \Big)^2         ds  \geq -    (L_{12}^{s}(\Omega)(t) )^2.         
    \end{equation} 
Therefore,  we have 
  $$
  L_{12}^{s}(\Omega)(t) \geq \frac{L^s_{12}(\Omega)(0)}{1+tL^s_{12}(\Omega)(0)},
  $$
which gives 
   $$
  \int_0^{t}L_{12}^{s}(\Omega)(\tau)  \geq   \log(1+tL^s_{12}(\Omega)(0)). 
  $$
 Hence, this complete the lower bound \eqref{L12_Low_est} gives the result.

\end{proof}

  \begin{proposition} \label{LOM_BMO} 
   Let $\Omega(t)$ be a solution to the leading order model \eqref{LOMeq} with   initial data \eqref{LOMeqOmega_0}. Then,   there exist constant $c$ independent of $\alpha$ such that we have the following estimates 
 \begin{equation} \label{LOM_BMO_Est}
|\Omega(t)|_{BMO}   \geq     c     \log(1+c\frac{t}{\alpha}),
\end{equation} 

\noindent on the time scale  $0\leq t \leq T(\alpha)=c \alpha \log(c |\log(\alpha)|) $.

   \end{proposition}
\begin{proof}
   Recall that      
   $$
   \Omega_0(r,\theta) = \frac{1}{\beta(\alpha)} r ^{\alpha}  \cos( \beta(\alpha) \log(r))+ r ^{\alpha} \sin(2\theta).
   $$
We proceed similarly to Proposition \ref{L12_Up_Low}. Setting $\Phi^{-1}(t,r,\theta)=(\Phi^{-1}_r (t,r,\theta),\Phi^{-1}_{\theta}(r,\theta)))=(r(t),\theta(t))$. From  Proposition \ref{L12_Up_Low}, we obtained the following by studying the characteristic of the leading order mode 
 \begin{equation}\label{FlowMapEqBMO}
 (r(t))^2 |\sin(2\theta(t))|   =     r^2  |\sin(2\theta)| .
\end{equation} 
Furthermore, to shorten  the notation as in Proposition \ref{L12_Up_Low}, we denote  $r(t)=re^{X(t)}$, where  
 \begin{equation}\label{XflowBMO}
  X(t,r,\theta)=2\int_0^tL_{12}(r(\tau))\cos(2\theta(\tau)) d\tau.
\end{equation} 
 We can now write the solution as 
      \begin{align}
   \Omega(t,r,\theta) &=\frac{1}{\beta(\alpha)}    (r(t))^{\alpha}  \cos(  \beta(\alpha) \log(r(t)))+ (r (t))^{\alpha} \sin(2\theta(t))\label{OmegaSol1}\\   
                                &=\frac{1}{\beta(\alpha)}  r^{\alpha}e^{\alpha  X(t)}  \Big( \cos( \beta(\alpha)   \log(r )) \cos(  \beta(\alpha)X(t) ) - \sin( \beta(\alpha)\log(r )) \sin( \beta(\alpha) X(t) )  \Big)   +   r^{\alpha}e^{\alpha X(t)} \sin(2\theta(t)) \label{OmegaSol2} 
 \end{align}
  \noindent   Using  \eqref{L12_Up_est} from Proposition  \ref{L12_Up_Low}, we have 
 \begin{align}\label{FlowUpBnd}
|X(t)|_{L^{\infty}} \leq c\frac{t}{\alpha}, 
 \end{align}

\noindent on the time $T(\alpha)=c \alpha \log(c |\log(\alpha)|)$, and since $\beta(\alpha)=|\log(|\log(\alpha)|^{1+\eta})|^{-1}$, we have   
 \begin{align}\label{ExpansionEst}
|\beta(\alpha) X(t)|, |\alpha X(t)|<1
 \end{align}
Thus, we can expand the following terms   
 \begin{align}
\cos(\beta(\alpha)X(t))) &= 1-(\beta(\alpha))^2 (X(t))^2 + \dots  \label{cosE} \\
\sin(\beta(\alpha)X(t)))&=\beta(\alpha)X(t) -  (\beta(\alpha))^3 (X(t))^3 + \dots  \label{sinE}\\
e^{\alpha X}&=1+ \alpha X  + \dots \label{expE}  
\end{align}
Hence, we can write 
 \begin{align}\label{OmegaExpanded}
   \Omega(t,r,\theta)    &=\frac{1}{\beta(\alpha)}  r^{\alpha}    \cos( \beta(\alpha)   \log(r ))   -r^{\alpha}  \sin( \beta(\alpha)\log(r ))   X(t)  +\beta(\alpha) E(r,\theta,\alpha)    +   r^{\alpha} \sin(2\theta(t)),
 \end{align}
where $E$  contains all the remaining terms from the expansion.  From  \eqref{ExpansionEst}, the term $E$ is uniformly bounded, and since the term   $r^{\alpha} \sin(2\theta(t))$ is also bounded, both terms can be absorbed when estimating $BMO$ norm. Thus, the main contribution comes from the first two terms in \eqref{OmegaExpanded}. Hence, we can start estimating the $BMO$ norm as follows:
 \begin{align}\label{OmegaBMOEst}
 | \Omega(t)|_{BMO}\geq \frac{1}{\pi r^2}\int_{B_r(0)} \Big|\Omega(t,s,\theta)-\frac{1}{\pi r^2}\int_{B_r(0)}\Omega(t,h,\sigma)  h dh d\sigma \Big| sds d\theta,
   \end{align}
  for any $r>0$. Now for the inner average integral, we observe from \eqref{XflowBMO} that due to symmetry, we have
     
   $$
    \frac{1}{\pi r^2}\int_{B_r(0)} h^{\alpha}  \sin( \beta(\alpha)\log(h ))   X(t,h,\sigma)   h dh d\sigma=0
   $$
  Furthermore, we have 
  
    $$
    \frac{1}{\pi r^2}\int_{B_r(0)}  \frac{1}{\beta(\alpha)}  h^{\alpha}    \cos( \beta(\alpha)   \log(h ))   h dh d\sigma= \frac{1}{\beta(\alpha)}  r^{\alpha}    \cos( \beta(\alpha)   \log(r )) +\bar{E}(r,\alpha),
   $$
   where $\bar{E}(r,\alpha)$ is uniformly bounded which can be absorbed. Therefore, we can estimate the $BMO$ norm as follows: 
   $$| \Omega(t)|_{BMO}\geq \frac{c}{\pi r^2}\int_{B_r(0)} \Big| \frac{1}{\beta(\alpha)}  s^{\alpha}    \cos( \beta(\alpha)   \log(s ))   - \sin( \beta(\alpha)\log(s ))   X(t) - \frac{1}{\beta(\alpha)}  r^{\alpha}    \cos( \beta(\alpha)   \log(r )) \Big| sds d\theta$$
  Now since the function  $\frac{1}{\beta(\alpha)}  r^{\alpha}    \cos( \beta(\alpha)   \log(r ))$ is purely radial, it is  $O(1)$ in $BMO$. Namely, from Lemma \ref{TypeII}, we have 
     
  $$\frac{c}{\pi r^2}\int_{B_r(0)} \Big| \frac{1}{\beta(\alpha)}  s^{\alpha}    \cos( \beta(\alpha)   \log(s ))   - \frac{1}{\beta(\alpha)}  r^{\alpha}    \cos( \beta(\alpha)   \log(r )) \Big| sds d\theta \leq c$$
  where $c$ is \textit{independent} of $\alpha$. Hence, we obtain 
   $$| \Omega(t)|_{BMO}\geq \frac{c}{\pi r^2}\int_{B_r(0)} \Big|  s^{\alpha}   \sin( \beta(\alpha)\log(s ))   X(t) \Big| sds d\theta$$ 
Considering the characteristic  of the inverse flow map
 $$
\dot{r}(t)=   2r(t) L^s_{12}(\Omega) \cos(2\theta(t)), \quad \dot{\theta}(t)=-2 L^s_{12}(\Omega)  \sin(2\theta(t)).$$
 We observe that dynamic of the trajectories  are hyperbolic from \eqref{FlowMapEqBMO}, and the angular sector $[0,\frac{\pi}{8}]$ is a trapping region. Thus, using \eqref{XflowBMO} gives
   \begin{align}\label{OmegaBMOEst}
| \Omega(t)|_{BMO}&\geq \frac{c}{\pi r^2}\int_0^{r} \int_0^{\frac{\pi}{8}} \Big|  s^{\alpha}   \sin( \beta(\alpha)\log(s ))    \int_0^tL_{12}(s(\tau))\cos(2\theta(\tau)) d\tau \Big| sds d\theta\\
&\geq \frac{c}{\pi r^2}\int_0^{r} \int_0^{\frac{\pi}{8}}   s^{\alpha}   |\sin( \beta(\alpha)\log(s ))|    \int_0^tL_{12}(s(\tau))  d\tau \Big| sds d\theta\\
&\geq \frac{c}{ r^2}\int_0^{r}   s^{\alpha+1}   |\sin( \beta(\alpha)\log(s ))|    \int_0^tL_{12}(s(\tau))  d\tau  sds,  
   \end{align}
 We obtained  from  \eqref{L12_Low_est} in Proposition \ref{L12_Up_Low} that 
$$\int_0^t L_{12}(\tau)    d\tau \geq c \log(1+c\frac{t}{\alpha}(1-s^{\alpha})).$$
Therefore, by choosing  $r=e^{-\frac{1}{\alpha}}$, we have 
  \begin{align}\label{OmegaBMOEst2}
| \Omega(t)|_{BMO}&\geq  
   c e^{\frac{2}{\alpha}}  \log(1+c\frac{t}{\alpha}) \int_0^{e^{-\frac{1}{\alpha}}}  s^{\alpha+1}   |\sin( \beta(\alpha)\log(s ))|     ds.  
   \end{align}
   Hence, it just remains to show that 
$$
   \int_0^{e^{-\frac{1}{\alpha}}}  s^{\alpha+1}   |\sin( \beta(\alpha)\log(s ))|     ds= c e^{-\frac{2}{\alpha}},
$$
 which we can  estimate in  the same manner as in Lemma \ref{BMOexLow}. Namely, we have 
  \begin{align}\label{OmegaBMOEst3}
 \int_0^{e^{-\frac{1}{\alpha}}}  s^{\alpha+1}   |\sin( \beta(\alpha)\log(s ))|     ds=\frac{1}{\beta(\alpha)}    \int_{\frac{\beta(\alpha)}{\alpha}}^{\infty}    e^{-\frac{u}{\beta(\alpha) } (\alpha+2)}  | \sin( u)  | du. 
  \end{align}
 \noindent  For any positive integer $k$, we observe
 
    \begin{align*} 
  \int_{\frac{3}{4}\pi+k\pi}^{\infty}  e^{-\frac{u}{\beta(\alpha) } (\alpha+2)}  | \sin(u)  |         du   & \geq   \sum_{j=k}^{\infty}  \int_{\frac{3}{4}\pi+j\pi}^{\frac{5}{4}\pi+j\pi}  c e^{-\frac{u}{\beta(\alpha) } (\alpha+2)}          du \\
&   \geq  c \frac{e^{-\frac{\frac{5}{4}\pi+k\pi}{\beta(\alpha) } (\alpha+2)}}{1-e^{-\frac{1}{\beta(\alpha) } (\alpha+2)}  } \geq c\beta(\alpha) e^{-\frac{\frac{5}{4}\pi+k\pi}{\beta(\alpha) } (\alpha+2)}
   \end{align*}
   where to obtain the last inequality, we used that    $\beta(\alpha)\rightarrow 0$ as $\alpha \rightarrow 0$. Thus, by appropriately choosing $\alpha$ , since $\frac{\beta(\alpha)}{\alpha} \rightarrow \infty$ as $\alpha \rightarrow 0$, we have    
   \begin{align}\label{OmegaBMOEst4}
 \int_{\frac{\beta(\alpha)}{\alpha}}^{\infty}    e^{-\frac{u}{\beta(\alpha) } (\alpha+2)}  | \sin( u)  | du \geq    c\beta(\alpha)  e^{-\frac{1}{ \alpha  } (\alpha+2)}.   
    \end{align}
     Hence, from \eqref{OmegaBMOEst2}, \eqref{OmegaBMOEst3}, and \eqref{OmegaBMOEst4}, we obtain the desired result       
    $$| \Omega(t)|_{BMO} \geq c    \log(1+ c  \frac{t}{\alpha}).    $$

   \end{proof}

\section{Leading Order Model Useful Estimates}\label{Sec_LOM_UE}
In this section, we prove some upper bound estimates on the leading order model. These will be used in Section \ref{Sec_BMO_remainder_Est} to estimate the error between the model and the 2d Euler equation. 
\begin{lemma} \label{LOMLcUpbnd}

   Let $\Omega(t)$ be a solution to the leading order model \eqref{LOMeq} with   initial data \eqref{LOMeqOmega_0}. Then,   there exist constants $c$ independent of $\alpha$ such that we have the following estimates 
    
\begin{equation} \label{L12cUpbnd}
| L_{12}^{c}(\Omega)(t)|_{L^{\infty}}  \leq    \frac{ct}{\alpha \beta(\alpha)}, \,\, \text{and}  \,\, \, | L_{12}^{c}(\Omega)(t)|_{C^{\alpha}}  \leq    \frac{ct}{\alpha  \beta(\alpha)}.
\end{equation}

   \end{lemma}
   \begin{remark} Similar type of estimate was obtain in \cite{SK} in the context of ill-posedness of 2d Euler in Besov regularity. 
   \end{remark}

 \begin{proof}
  \noindent Here we would like to estimate the following:  
 
   $$
  L_{12}^{c}(\Omega) =  c\int_{r}^{1}\int_{0}^{\frac{\pi}{2}}\frac{  \Omega(t,r,\theta)  }{s}     \cos(2\theta) d\theta ds.  
  $$
 We begin the argument  similar to    Proposition \ref{LOM_BMO}.  From \eqref{OmegaSol2}, we can write the solution as: 
 \begin{align}
   \Omega(t,r,\theta)    &=\frac{1}{\beta(\alpha)}  r^{\alpha}e^{\alpha  X(t)}  \Big( \cos( \beta(\alpha)   \log(r )) \cos(  \beta(\alpha)X(t) ) - \sin( \beta(\alpha)\log(r )) \sin( \beta(\alpha) X(t) )  \Big)   +   r^{\alpha}e^{\alpha X(t)} \sin(2\theta(t))
 \end{align}
 where recall that $X(t,r,\theta)=2\int_0^tL_{12}(r(\tau))\cos(2\theta(\tau)) d\tau $. Again, similar to    Proposition \ref{LOM_BMO},  using \eqref{FlowUpBnd}, $|X(t)|_{L^{\infty}} \leq c t \alpha^{-1}$, on our time scale $T(\alpha)$, we have $|\beta(\alpha) X(t)|, |\alpha X(t)|<1$. Thus, we can expand our solution using \eqref{cosE},  \eqref{cosE}, and \eqref{expE}, as follows:
   \begin{align}\label{OmegaExpanded2}
   \Omega(t,r,\theta)    &=\frac{1}{\beta(\alpha)}  r^{\alpha}    \cos( \beta(\alpha)   \log(r ))   -r^{\alpha}  \sin( \beta(\alpha)\log(r ))   X(t)  +\beta(\alpha) E(r,\theta,\alpha)    +   r^{\alpha} \sin(2\theta(t)),
 \end{align}
 where $E$ is a  bounded error term.  When computing $L_{12}^{c}(\Omega)$, by symmetry, we observe that the terms $(\beta(\alpha))^{-1}  r^{\alpha}    \cos( \beta(\alpha)   \log(r )) $, and $r^{\alpha} \sin(2\theta(t))$ vanish. In addition,  since the error term $E$ is bounded and contains a  factor of  $\beta(\alpha)$, it will be absorbed. Therefore, it suffices to estimate the following
    \begin{align}\label{L12cOsc_}
   |  L_{12}^{c}(\Omega)  |  \leq \Big| c\int_{r}^{1}\int_{0}^{\frac{\pi}{2}} \frac{s^{\alpha}}{s}    \sin( \beta(\alpha)\log(r ))   X(t)     \cos(2\theta) d\theta ds  \Big|. 
  \end{align}
  
 \noindent We would like to recast the above integral  into an oscillatory  integral so we can gain a better estimate. From the definition of $X(t)$, we set $\sigma(r,\tau)= \int_{0}^{\frac{\pi}{2}}  \cos(2\theta(\tau)) \cos(2\theta) d\theta   $, which is  bounded, and rewrite \eqref{L12cOsc_} as follows: 
 \begin{align}\label{L12cOsc0}
  |   L_{12}^{c}(\Omega)  |   \leq \Big| c \int_0^t \int_{r}^{1} \frac{s^{\alpha}}{s}      \sin( \beta(\alpha)\log(r ) )L_{12}(r(\tau))  \sigma(r,\tau) ds  d\tau \Big|. 
\end{align}
  
 \noindent First we focus on the radial integral, so we set $\lambda=\beta(\alpha) |\log(r)|$ and   introduce the change of variable   $u= \log(s) (\log(r))^{-1}$ that gives 
 \begin{align}\label{L12cOsc1}
   \int_{r}^{1} \frac{s^{\alpha}}{s}     \sin( \beta(\alpha)\log(r ))  2L^{s}_{12}(s(\tau))   \sigma(\tau,s)   ds= c \frac{\lambda}{\beta(\alpha)}  \int_{0}^{1}e^{-\frac{\alpha}{\beta(\alpha)} \lambda u}    \cos( \lambda u)    L^s_{12}(u(\tau))   \sigma(\tau , u  )  du . 
\end{align}

 \noindent  Since $\sigma$ is bounded, and $  \alpha L_{12}(\Omega)$  is also bounded from   	\eqref{L12_Up_est} in Proposition \ref{L12_Up_Low}, we can obtain the following by observing that $\lambda= \beta(\alpha)  |\log(r)| \rightarrow \infty$ as     $r \rightarrow 0$:
$$
  \bigg|  \int_{0}^{1}e^{-\frac{\alpha}{\beta(\alpha)} \lambda u}    \cos( \lambda u)    \alpha L^s_{12}(u(\tau))   \sigma(\tau , u  )  du      \bigg|    \leq \frac{c}{\lambda}.
  $$
 Thus, from  $\eqref{L12cOsc0}$ and  $\eqref{L12cOsc1}$, we have  
  
  $$
  | L_{12}^{c}(\Omega)|_{L^{\infty}}   \leq  \frac{c t}{\alpha \beta(\alpha)},
  $$
  
\noindent and  the $C^{\alpha}$ estimate follows in a similar manner, and therefore this completes the proof.

  \end{proof}

  \begin{lemma} \label{Omega_Grad_HBnd}

Let $\Omega(t)$ be a solution to the leading order model \eqref{LOMeq} with   initial data \eqref{LOMeqOmega_0}. Then,  we have the following estimate:  
\begin{equation}\label{Omega_H_est}  
|  \Omega(t)|_{C^{\alpha}}   \leq \frac{1}{\beta(\alpha)}   e^{\frac{ c t}{\alpha}}, \quad  | \dth \Omega(t)|_{C^{\alpha}}  \leq  c    e^{\frac{ct}{\alpha}},  \quad  | r \dr \Omega(t)|_{C^{\alpha}}  \leq     c    e^{\frac{ct}{\alpha}},
\end{equation} 
and 
 \begin{equation} \label{OmegaCirc4BMO}
 |\Omega|_{\mathring{C}^{1,\alpha}}   \leq       \frac{c}{\beta(\alpha)}   e^{\frac{ c t}{\alpha}},
\end{equation} 
where   $c$ is independent of $\alpha$.

   \end{lemma}
   
   \begin{proof}
   Recall that our initial data is:
      $$
   \Omega_0(r,\theta) = \frac{1}{\beta(\alpha)} r ^{\alpha}  \cos( \beta(\alpha) \log(r))+ r ^{\alpha} \sin(2\theta).
   $$
 Thus, we have the following  
    $$
|  \Omega_0|_{C^{\alpha}} \leq \frac{1}{\beta(\alpha)}   , \quad    |\dth  \Omega_{0}|_{C^{\alpha}} \leq c , \quad \text{and}   \quad   |r \dr \Omega_{0}|_{C^{\alpha}} \leq c  ,
   $$  
   where $c$ is independent of $\alpha$. Here we used that $\alpha (\beta(\alpha))^{-1} \leq c$ for all $\alpha$.  Now because the vorticity is being transported, we have 
   $$
   \Omega(r,t,\theta)=\Omega_{0}(\Phi^{-1}_t(r,\theta)).
   $$
  Hence, from \eqref{L12_Up_est} in Proposition \ref{L12_Up_Low}, and  standard transport estimates, we can obtain the following  
\begin{equation}\label{Omega_H_est_}  
|  \Omega(t)|_{C^{\alpha}}   \leq \frac{c}{\beta(\alpha)}   e^{ \frac{ct}{\alpha}}, \quad  | \dth \Omega(t)|_{C^{\alpha}}  \leq  c    e^{\frac{ct}{\alpha}}, \,\,\, \text{and}  \,\,\, | r \dr \Omega(t)|_{C^{\alpha}}  \leq     c    e^{\frac{ct}{\alpha}}.
\end{equation} 
with  $c$ is independent of $\alpha$. Recall from subsection \ref{SubSec_Notation} that $\mathring{C}^{1,\alpha}$ norm is defined as follows:
   \begin{equation}  
   |\Omega|_{\mathring{C}^{1,\alpha}}=    |\Omega|_{C^{\alpha}}+   |\dth  \Omega|_{C^{\alpha}}+|r \dr \Omega|_{C^{\alpha}}. 
\end{equation} 
Thus, the  $|\Omega|_{\mathring{C}^{1,\alpha}}$ estimates follows from \eqref{Omega_H_est_}, and hence this completes the proof. 
 
   \end{proof}

\section{Remainder Estimate } \label{Sec_BMO_remainder_Est}

In this section, we estimate the difference  between the Euler equation solution $\omega(t)$ and the leading order model solution $\Omega(t)$. Recall that the full 2d Euler, $\omega(t)$, solves:
$$\dt \omega+\frac{1}{r} \dr \psi \dth \omega-\frac{1}{r}\dth \psi \dr\omega=0$$ 
Now the solution to the Euler equation can be decomposed as follows:
 \begin{equation}\label{EulerDecompose} 
\omega(t)=\omega_r(t)+ \Omega(t),
\end{equation}  
where $\omega_r(t)$, we call  remainder, is the difference  between the Euler equation, $\omega(t)$,  and 
 the leading order model, $\Omega(t)$ which is:
 \begin{equation}\label{LMO_RemainderS}  
 \dt \Omega+  2 L^s_{12}(\Omega) \sin(2\theta)  \dth \Omega- 2 rL^s_{12}(\Omega) \cos(2\theta) \dr \Omega =0.
\end{equation} 
The stream function can be written as follows:  $$\psi(\omega)=\psi(\omega_r)+ \psi(\Omega),$$ 
 From the Biot-Savart law decomposition (Proposition \ref{ElgindiDecomp}), and Remark \ref{Remark_psi_remainder}, we can write:

\begin{equation} \psi_r(\Omega)=\psi(\Omega)-r^2L^{s}_{12}(\Omega)(r) \sin(2\theta)-r^2L^{c}_{12}(\Omega)(r) \cos(2\theta),
 \end{equation} 
Therefore, we obtain the following equation for the remainder $\omega_r(t)=\omega(t)- \Omega(t) $

  \begin{equation}\label{RemainderEqBMO}
\begin{split}
  \dt \omega_r& -\Big( 2rL^s_{12}(\Omega)(r)\cos(2\theta)- 2rL^c_{12}(\Omega)(r)\sin(2\theta)+ \frac{1}{r}   \dth \psi_r(\Omega) + \frac{1}{r}   \dth \psi(\omega_r) \Big) \dr \omega_r 
   \\ &+ \Big(2L^s_{12}(\Omega)(r)\sin(2\theta)+ 2L^c_{12}(\Omega)(r)\cos(2\theta)   + \frac{1}{r}    \dr \psi_r(\Omega)+ \frac{1}{r}   \dr \psi(\omega_r)\Big) \dth \omega_r \\
  &-    \Big( \frac{1}{r}    \dth \psi(\omega_r)  - 2rL^c_{12}(\Omega)(r)\sin(2\theta)+ \frac{1}{r}   \dth \psi_r(\Omega)   \Big)   \dr \Omega\\ 
  &+  \Big(  \frac{1}{r}    \dr \psi(\omega_r)   + 2rL^c_{12}(\Omega)(r)\cos(2\theta)+   \frac{1}{r}    \dr \psi_r(\Omega) \Big)   \dth \Omega=0.
\end{split}
\end{equation}
with $\omega_r|_{t=0}=0$.

\begin{proposition}\label{RemainderEqBMOP}
Let $\omega_r(t)$ be  the solution to \eqref{RemainderEqBMO} with $\omega_r|_{t=0}=0$, then for some $0<\delta<1$, we have the following:
\begin{equation}\label{RemainderEqestBMO}
|\omega_r(t)|_{C^{\alpha}} \leq c \alpha^{\delta}, 
\end{equation}
on the time interval $0 \leq t\leq  T(\alpha)=c \alpha \log(c |\log(\alpha)|) $ with $c$ independent of $\alpha$. 

\end{proposition}

\begin{proof} Given that  $\omega_r$ satisfies \eqref{RemainderEqBMO}, we compose with the flow map of the transport  estimate the  H\"older norm as follows:
 \begin{equation}\label{RemainderEq_H_BMO} 
\begin{split}
  \frac{d}{dt}|\omega_r|_{C^{\alpha}}  &\leq c \Big( |\frac{1}{r^2}    \dth \psi(\omega_r) |_{C^{\alpha}} + |L^c_{12}(\Omega)|_{C^{\alpha}} +   | \frac{1}{r^2}   \dth \psi_r(\Omega)  |_{C^{\alpha}}   \Big)  |r \dr \Omega|_{C^{\alpha}} \\ 
 &+c  \Big(  | \frac{1}{r} \dr \psi(\omega_r)|_{C^{\alpha}}+ |L^c_{12}(\Omega) |_{C^{\alpha}}  +   | \frac{1}{r} \dr \psi_r(\Omega)|_{C^{\alpha}}     \Big)   |\dth \Omega|_{C^{\alpha}} 
\end{split}
\end{equation}
When we do the Schauder estimates for the above equation,  we can drop $L^{1}$ term since the radius of support of our solution will be bounded on our time scale. We will estimate each similar group of terms together: 

\noindent   \textbf{Terms  $| \frac{1}{r^2}  \dth \psi(\omega_r) |_{C^{\alpha} }$ and $|\frac{1}{r} \dr \psi(\omega_r)|_{C^{\alpha}}$} 

\noindent We obtain the estimate on these terms  by  standard elliptic estimates 
\begin{align}
& |\frac{1}{r^2} \dth \psi(\omega_r) |_{C^{\alpha}} \leq \frac{c}{\alpha}|\omega_r|_{C^{\alpha}}  . \label{dthPsRemiEst_BMO}\\
& |\frac{1}{r} \dr \psi(\omega_r) |_{C^{\alpha}} \leq \frac{c}{\alpha}|\omega_r|_{C^{\alpha}}  . \label{drPsRemiEst_BMO}
\end{align}

\noindent  \textbf{Term  $ |L^c_{12}(\Omega)|_{C^{\alpha}}$   } 
 
\noindent  This term can be esimated by using Proposition \ref{LOMLcUpbnd}.    
\begin{align}\label{L12cRemEstBMO}
& |L^c_{12}(\Omega)|_{C^{\alpha}}  \leq    c e^{\frac{t}{\alpha \beta(\alpha)}}.   
\end{align}

\noindent \textbf{Terms   $|\frac{1}{r} \dr \psi_r(\Omega)|_{C^{\alpha}}$ and $| \frac{1}{r^2}  \dth \psi_r(\Omega) |_{C^{\alpha}}$  } 

\noindent We estimate  these terms by using Proposition \ref{ElgindiDecomp}, and  \eqref{OmegaCirc4BMO} from Lemma  \ref{Omega_Grad_HBnd} 
\begin{align}
&|\frac{1}{r}\dr \psi_r(\Omega)|_{C^{\alpha}}\leq c    |\Omega|_{\mathring{C}^{1,\alpha}}  \leq   \frac{c}{\beta(\alpha)}   e^{\frac{ c t}{\alpha}} \label{drPsirRemEst_BMO},  \\
& |\frac{1}{r^2} \dth \psi_r(\Omega) |_{C^{\alpha}} \leq    c |\Omega|_{\mathring{C}^{1,\alpha}}  \leq       \frac{c}{\beta(\alpha)}   e^{\frac{ c t}{\alpha}}.  \label{dthPsirRemEst_BMO}
\end{align}

\noindent  \textbf{Terms $|\dth \Omega|_{C^{\alpha}}$ and $|r \dr \Omega|_{C^{\alpha}}$  } 

\noindent  These terms we estimate them by using \eqref{Omega_H_est} from  Lemma \ref{Omega_Grad_HBnd}.  
\begin{align}
&|\dth \Omega|_{C^{\alpha}} \leq  c   e^{\frac{ct}{\alpha}} ,\label{dthOmegaRemEst_BMO} \\
&|r \dr \Omega|_{C^{\alpha}} \leq   c   e^{\frac{ct}{\alpha}} .\label{drOmegaRemEst_BMO} 
 \end{align}

\noindent   \textbf{Total Remainder Estimate}

\noindent Using  the remainder equation \eqref{RemainderEq_H_BMO} with \eqref{dthPsRemiEst_BMO}, \eqref{drPsRemiEst_BMO},  \eqref{L12cRemEstBMO}, \eqref{drPsirRemEst_BMO}, \eqref{dthPsirRemEst_BMO}, \eqref{dthOmegaRemEst_BMO}, and \eqref{drOmegaRemEst_BMO}, we have 
 \begin{align*}
  \frac{d}{dt}|\omega_r(t)|_{C^{\alpha}} & \leq  \Big(  \frac{c}{\alpha}|\omega_r(t)|_{C^{\alpha}} +   \frac{c}{\beta(\alpha)}   e^{\frac{ c t}{\alpha}}+  c e^{\frac{t}{\alpha \beta(\alpha)}}     \Big)   c    e^{\frac{ct}{\alpha}}.  
  \end{align*}
Therefore, by using Gronwall inequality with the fact that   $\omega_r|_{t=0}=0$, we have   
\begin{align*}
   |\omega_r(t)|_{C^{\alpha}} & \leq \Big( \int_0^t         (\frac{c}{\beta(\alpha)}   e^{\frac{ c \tau}{\alpha}}+  c e^{\frac{\tau}{\alpha \beta(\alpha)}} )d\tau \Big) \exp{\Big( \int_0^t   \frac{c}{\alpha}  e^{\frac{\tau}{\alpha}}   } d\tau\Big)  \leq   c    \Big(          \frac{ \alpha}{\beta(\alpha)}   e^{\frac{ c t}{\alpha}}+   \alpha \beta(\alpha) e^{\frac{c t}{\alpha \beta(\alpha)}}   \Big)   \exp{\Big(    c   e^{\frac{t}{\alpha}}   }\Big) . 
  \end{align*}
  To close the estimate, we choose  the time scale  $T(\alpha)= c \alpha \log(c |\log(\alpha)|) $ for $c$ small enough and independent of $\alpha$, and since $\beta(\alpha)=  \log( |\log(\alpha)|^{1+\eta}) $, we have 
$$
   |\omega_r(t)|_{C^{\alpha}}     \leq  c  \alpha^{\delta},  
$$
for $0\leq t \leq  T(\alpha)$,  and some $0<\delta<1$.  This gives the result.

   \end{proof}

   \section{Main Result  }\label{Sec_Main}
In this section, we prove our main result. We recall the statement of main result for convenience of the reader. 
\begin{theorem*} (Strong ill-posedness for vorticity in $BMO$):  For any $\epsilon>0$, and $0<\alpha<1$, there exist compactly supported initial vorticity  $\omega_0^{\alpha,\epsilon} $ on $\mathbb{R}^2$  such that the unique solution to the 2d Euler equation  satisfies the following: 
$$|\omega_0|_{BMO}=\epsilon, \,\,  \text{and} \,\,  \sup_{0\leq t \leq T(\alpha)}|\omega(t)|_{BMO} \geq   \log(c\log(c |\log(\alpha)|)), $$ 
  where 
  $T(\alpha)=c \alpha \log(c |\log(\alpha)|) $ and c is independent of $\alpha$.    
  \end{theorem*}

  \begin{proof}
  
  Recall   from  \eqref{LOMeqOmega_0} that  the initial vorticity is
  \begin{equation}\label{MT_Omega_0}
  \omega_0(r,\theta)=\Omega_0(r,\theta) =\frac{1}{\beta(\alpha)}r^{\alpha} \cos(\beta(\alpha)\log(r))+  r ^{\alpha} \sin(2\theta)
\end{equation} 
compactly supported on the ball of radius 1 and with  $\beta(\alpha)= |\log(|\log(\alpha)|^{1+\eta})|^{-1}$ for some $\eta>0$. In addition, recall from  \eqref{EulerDecompose} that we can write the solution of the 2d  Euler equation as:
   \begin{equation}\label{EulerSolution}
  \omega(t)=\Omega(t)+\omega_r(t)
  \end{equation}
  where $\Omega(t)$ is the solution to the leading order model \eqref{LMO_RemainderS} with initial data $\Omega|_{t=0}= \omega_0(r,\theta)$ in \eqref{MT_Omega_0}, and $\omega_r(t)$ is the remainder  solving \eqref{RemainderEqBMO}   with initial data   $\omega_r|_{t=0}=0$. From Proposition \ref{RemainderEqBMOP}, we have $|\omega_r(t)|_{C^{\alpha}} \leq c \alpha^{\delta}$, on the time scale $T(\alpha)=c \alpha \log(c |\log(\alpha)|)$, which implies  that  $|\omega_r(t)|_{L^{\infty}} \leq c \alpha^{\delta}$, and hence we have 
   \begin{equation}\label{RemainderSolutionBMOEst}
  |\omega_r(t)|_{BMO} \leq c,
  \end{equation}
  for a constant $c$  independent of $\alpha$. Further, from Proposition \eqref{LOM_BMO},  we have 
   \begin{equation}\label{LOMBMOSolutionEst}
|  \Omega(t)  |_{BMO} \geq  c          \log(1+c\frac{t}{\alpha}), 
  \end{equation}
on our time interval $0\leq t \leq T(\alpha)=c\alpha \log(c |\log(\alpha)|) $. Thus, from \eqref{EulerSolution}, \eqref{RemainderSolutionBMOEst}, and \eqref{LOMBMOSolutionEst} we obtain that the vorticity of the full 2d  Euler equation satisfies 
 $$\sup_{0\leq t \leq T(\alpha)}| \omega (t)|_{BMO} \geq  c \log(c\log(c |\log(\alpha)|)),$$
   where $T(\alpha)=c\alpha \log(c |\log(\alpha)|)   \rightarrow 0$ as $\alpha \rightarrow 0$. Now since the constants  are independent of $\alpha$, we can  rescale the initial data by $\epsilon>0$, and obtain the strong ill-posedness result.
 
\end{proof}

\section*{Acknowledgements}
I would like to thank Professor Tarek M. Elgindi for important discussions and feedback. I also would  like to thank Professor Theodore D. Drivas and Professor John Anderson for their valuable feedback.

\vskip 0.3 cm

Mathematics Department, Stony Brook University, Stony Brook NY, 11794
 
\textit{Email address:} karim.shikhkhalil@stonybrook.edu

\end{document}